\documentclass[11pt]{article}
\usepackage{amssymb}
\usepackage{amsmath}
\usepackage{amsthm}
\usepackage{graphicx}
\usepackage{graphics}
\usepackage{amsthm}

\input pictex.tex

\newcommand{\s}{\sigma}
\newcommand{\vsp}{\vspace{0.1cm}}
\theoremstyle{definition}

\newtheorem{thm}{Theorem}[section]

\newtheorem{prop}[thm]{Proposition}

\newtheorem{lem}[thm]{Lemma}

\newtheorem{rem}[thm]{Remark}

\newtheorem{ex}[thm]{Example}

\input epsf

\begin{document}

\date{}
\author{Andr\'es Navas}

\title{Groups, Orders, and Laws}
\maketitle


\vspace{-0.5cm}



\section{Introduction}

In the general program of understanding the group of diffeomorphisms of a given (compact) manifold 
as an infinite-dimensional analogue of a (simple) Lie group, several questions arise naturally, as for example:

\vspace{0.18cm}

\noindent {\em -- The Tits alternative:} Do all non-virtually-solvable subgroups contain (non-Abelian) 
free subgroups in two generators ?

\vspace{0.18cm}

\noindent {\em -- The Milnor alternative:} Does every finitely-generated subgroup 
have either polynomial or exponential growth ? In case of exponential growth, is it 
necessarily uniform ?

\vspace{0.18cm}

A complete answer to all of this (in arbitrary regularity and any dimension) is far from being available. 
Nevertheless, several relevant results are known, notably in the 1-dimensional case. Indeed, the 
growth alternative fails for groups of $C^1$ diffeomorphisms of 1-manifolds but holds in class 
$C^{1+\alpha}$ (yet the question about uniform growth remains open; see \cite{growth}). Concerning the 
Tits alternative, it is known to be false in class $C^{\infty}$ (a counter-example is given by the smooth 
realization of Thompson's group $F$ by diffeomorphisms of the interval \cite{GS}), though it is an 
open question in the real-analytic case. 

The validity of the Tits alternative obviously implies that if the underlying group satisfies a 
(nontrivial) {\em law}, then it must be virtually-solvable. (Notice, however, that for linear groups, this 
can be easily established without using Tits' theorem \cite{tits} just passing to the Zariski closure and using 
the classification of algebraic Lie groups.) In this work, we will consider this baby form of the Tits alternative 
for groups acting on 1-manifolds. More precisely, we deal with the next

\vspace{0.4cm} 

\noindent{\bf Question (i).} Let $\Gamma$ be a subgroup of the group $\mathrm{Homeo}_+(\mathbb{R})$
of orientation-preserving homeomorphisms of the real line. If $\Gamma$ satisfies a law, must 
it be virtually-solvable ?

\vspace{0.4cm}

This question can be addresed more generally for groups of homeomorphisms (diffeomorphisms) of 
any manifold. However, as every group acting on the interval can be realized as a group acting on any 
manifold --this can be made even preserving any $C^k$ regularity\footnote{Roughly, this is done as follows: 
we view a punctured disk $D$ in the manifold  as the union of a disjoint family of open intervals all arising 
from the origin; then on each interval we copy the original action, and we extend all the elements 
as being the identity outside $D$; to preserve regularity, prior to this we apply the Muller-Tsuboi 
trick \cite[Exercise 5.1.14]{book} in order to make all group elements flat at the endpoints.}--, 
the 1-dimensional case appears as the first nontrivial case to deal with. Notice moreover that 
in dimension 1, the case of (groups of homeomorphisms of) the circle reduces to that of the 
real line (equivalently, to the interval) due to the validity of the Ghys-Margulis weak form of 
the Tits alternative \cite{margulis,ghys} that we state as follows: If $\Gamma$ is a subgroup 
of $\mathrm{Homeo}_+(\mathrm{S}^1)$ without free subgroups, then its commutator 
subgroup acts with global fixed points. (Thus, after opening the circle at one of these 
fixed points, we obtain an action of $[\Gamma, \Gamma]$ on the real line...)

We strongly suspect a negative answer to the question above (perhaps even in the framework of 
amenable groups). Nevertheless, our goal is mostly to discuss the problem from several viewpoints, 
to provide partial positive results and eventual applications that are interesting by themselves, and 
to address many related questions. 

Every subgroup of $\mathrm{Homeo}_+(\mathbb{R})$ is obviously torsion-free. We hence start 
by noticing that it is already nontrivial to give examples of torsion-free groups satisfying 
laws and which are non-virtually-solvable. However, following a suggestion 
of Y. de Cornulier, one may produce a lot of examples by using the next classical 
lemma of Higman \cite{higman}: If $F$ is a free group and $N$ a normal subgroup, then 
$F / [N,N]$ is torsion-free. As a consequence, if $F/N$ satisfies a law $W (a_1,\ldots,a_n) = id$ 
and is non-virtually-solvable, then $F/[N,N]$ is torsion-free, non-virtually-solvable, and satisfies 
the law $\big[ W(a_1,\ldots,a_n),W(a_{n+1},\ldots,W(a_{2n})) \big] \!=\! id$. Recall, moreover, 
that there are many examples of non-virtually-solvable groups satisfying laws. Actually, there 
well-known non-amenable examples, as for instance the Burnside groups $B(n)$, with 
$n > 666$ an odd integer, defined by (see \cite{adyan})
$$B(n) := \big\langle a,b \!: W^n = id 
\mbox{ for every (nonempty) word } W \mbox{ in } a \mbox{ and } b \big\rangle.$$

\noindent{\bf Question (ii).} If we write $B(n) \sim F_2 / N$, does the group $F_2 / [N,N]$ 
embed into $\mathrm{Homeo}_+ (\mathbb{R})$ ?

\vspace{0.38cm}

There are many other sources of motivation for Question (i) above, most of which come  from 
the theory of left-orderable groups. Recall that a group is said to be {\em left-orderable} if it 
admits a total order relation that is invariant under left multiplication (a {\em left-order}, for 
short).  Such a group is necessarily torsion-free, but the converse is far from being true. Actually, 
it is a well-known fact that a countable group is left-orderable if and only if it embeds into 
$\mathrm{Homeo}_+(\mathbb{R})$; see \cite[Th. 6.8]{ghys-survey}. The next question was 
addressed by Linnell.

\vspace{0.4cm} 

\noindent{\bf Question \cite{linnell}.} Let $\Gamma$ be a left-orderable group. Assume $\Gamma$ 
contains no free subgroup. Is $\Gamma$ {\em locally indicable} ? (i.e. does every finitely-generated 
subgroup homomorphically surjects into $\mathbb{Z}$ ?)

\vspace{0.4cm}

Although the answer to this question is unclear, a major evidence pointing in the affirmative 
direction is a beatiful theorem of Morris-Witte \cite{morris}: 
Every amenable, left-orderable group is locally indicable. 
In our approach, we should stress that even for left-orderable groups satisfying a law, local 
indicability remains an open question. A relevant case concerns Engel type laws. To properly 
state this, given two elements $f,g$ in a group $\Gamma$, we let the $n^{\mathrm{th}}$-commutator 
$[f,g]_n$ of $f$ and $g$ be inductively defined by $[f,g]_1 := fgf^{-1}g^{-1}$ and 
$[f,g]_{k+1} := [[f,g]_{k},g]_1$. One says that the group satisfies the $n^{\mathrm{th}}$-Engel 
condition (or it is an $n$-Engel group) if $[f,g]_n = id$ for all $f,g$ in $\Gamma$, 
and that it is an Engel group if it is $n$-Engel for some $n \!\in\! \mathbb{N}$. 

Engel's condition remains somewhat mysterious beyond the framework of finite or linear groups. 
A central question of Plotkin asks whether every finitely-generated torsion-free Engel group 
is locally nilpotent. This has been shown to be true when $n \leq 4$ (even without 
the assumption of torsion-freenes; see \cite{havas}), but remains open in general. (Yet Juhasz 
and Rips have recently announced a negative answer for large-enough $n$.) Plotkin's 
question has been also addressed under stronger conditions than torsion-freeness. 
For example, a theorem of Kim and Rhemthulla establishes that bi-orderable Engel 
groups are nilpotent \cite{KR}. This has been extended under the weaker condition 
of local indicability \cite{glass}, but remains open for general left-orderable groups 
for $n > 4$ (see \cite{LM} for the case $n=4$). 

Locally indicable groups are particularly important in the theory of left-orderable groups. According 
to Conrad \cite{conrad} and Brodski \cite{brodski}, these groups are those that support a left-order 
with a property slightly weaker than bi-invariance (see \cite{fourier} for a simpler proof Brodski's theorem). 
Although the original Conrad's definition of this property is quite involved (see \cite{BR,KK}), it was noticed 
in \cite{leslie,fourier} that it is equivalent to a much simpler one: For every $f \succ id$ and $g \succ id$, 
one has $f g^2 \succ g$. This naturally suggests introducing the notion of {\em verbal property} 
for an order: Given a (nonempty) reduced word $W = W(a,b)$ carrying positive and negative exponents, 
we say that a left-order on a group is a {\em $W$-order} if $W(f,g)$ is always a positive element 
(i.e. larger than the identity) whenever both $f$ and $g$ are positive. For example, the $W$-orders 
for $W(a,b) := a^{-1} b a$ are the bi-invariant ones, and the Conradian orders are the $W$-orders 
for $W(a,b) := a^{-1} b a^2$. 

\vspace{0.4cm}

\noindent {\bf Question (iii).} Does there exist a word $W$ such that the $W$-orders are those that 
satisfy an specific and relevant algebraic property different from bi-orderability or local indicability ?

\vspace{0.4cm}

Conradian orders are special not only from a verbal viewpoint, but also from a dynamical one. Indeed, 
it has been shown in \cite{fourier}, and later pursued in \cite{new-order}, that the Conradian orders are 
those for which there is no resilient pair; more precisely, there are no group elements $f,g,h_1,h_2$ 
satisfying 
\begin{equation}\label{resilient}
h_1 \prec f^n h_1 \prec f^n h_2 \prec g^n h_1 \prec g^n h_2 \prec h_2
\end{equation}
for all $n \in \mathbb{N}$. (Notice that it suffices to check this for $n=1$ to garantee it for all 
$n \geq 1$, but this way of stating is more natural in view of the discussion below.) 
As a consequence, if a group admits non Conradian orders, then it must contain free subsemigroups;  
in particular, it cannot satisfy a semigroup law. (The latter facts were first established in \cite{LMR} 
using quite different arguments.) 

In our approach, what is relevant with condition (\ref{resilient}) is that it 
is somewhat related to our strategy to rule out laws in groups having a 
rich-enough action. Indeed, assume that for all $k \in \mathbb{N}$,  
there are group elements $f,g$ and $h_i,h_i',\bar{h}_i,\bar{h}_i'$, $i \in \{1,\ldots,2k+1\}$, 
such that for each nonzero integer $n$, 
\begin{equation}\label{n-resilient-1}
\mbox{either } \quad  
\bar{h}_{i-1} \prec f^n (h_i)  \prec f^n (h_{i}') \prec \bar{h}_{i-1}' 
\quad \mbox{ or } \quad  
\bar{h}_{i} \prec f^n (h_i)  \prec f^n (h_{i}') \prec \bar{h}_i' 
\end{equation}
whenever $i \!\in\! \{2,\ldots,2k+1\}$, and 
\begin{equation}\label{n-resilient-2}
\mbox{either } \quad  
h_{i} \prec g^n (\bar{h}_i)  \prec g^n (\bar{h}_{i}') \prec h_{i}' 
\quad \mbox{ or } \quad  
h_{i+1} \prec g^n (\bar{h}_i)  \prec g^n (\bar{h}_{i}') \prec h_{i+1}' 
\end{equation}
whenever $i \!\in\! \{1,\ldots,2k \}$. Then, as we will see below (c.f. Example 
\ref{law-order}), the underlying group cannot satisfy any law. This motivates still another 

\vspace{0.4cm}

\noindent {\bf Question (iv).} Does there exist an algebraic description of the set of orders for which 
(\ref{n-resilient-1}) and (\ref{n-resilient-2}) cannot hold simultaneously for any group elements 
and some $n > 2$ ? In particular, does it coincide with the set of $W$-orders for a certain 
word $W$ ? 


\section{Some results}

\noindent{\bf (a) Groups of diffeomorphisms satisfying laws.} As already mentioned, the 
Tits alternative fails in $\mathrm{Diff}^{\infty}_+([0,1])$, and it is an open question in 
the group $\mathrm{Diff}^{\omega}_+([0,1])$ of real-analytic diffeomorphisms. 
However, as we next state, its baby form remains true in the latter context. 

\vspace{0.4cm}

\noindent{\bf Theorem A.} {\em Every subgroup of $\mathrm{Diff}^{\omega}_+([0,1])$ 
satisfying a law is solvable.}

\vspace{0.4cm}

Notice that the result above forces virtually-solvable subgroups of $\mathrm{Diff}_+^{\omega}([0,1])$ to 
be solvable. (This was already known; see also \cite[Cor. 2]{bleak-universal} and \cite[Rem. 6.4]{navas-brasil}.)

It is very likely that Theorem A is still true for subgroups of 
$\mathrm{Diff}^{1+\alpha}_+([0,1])$. (The $C^1$ case is unclear to the author.) 
However, besides of a much more technical work, this would certainly require 
a different viewpoint, hence we will address this seemingly hard issue elsewhere. 
Just to give a flavour on the difficulties here, let us mention that Theorem A 
is still true for groups of piecewise-linear 
homeomorphisms of the interval. Indeed, if such a group is non-solvable, then a result of Bleak 
\cite{bleak-universal} establishes that it contains a copy of the group $G$ defined by 
$$G:= \bigoplus_{n \geq 0} G_n, \quad \mbox{where } G_0 := \mathbb{Z} 
\quad \mbox{and} \quad  G_{n+1}:= G_n \wr \mathbb{Z}.$$ 
Moreover, by a result of Akhmedov (see \cite[Lemma 2.1]{akhmedov-girth}), $G$ satisfies 
no nontrivial law. Unfortunately, our method of proof does not yield the fact 
that $G$ satisfies no law, despite the fact that $G$ can be realized as a group of 
$C^{\infty}$ diffeomorphisms of the interval (see \cite{navas-mexico}).   

\vspace{0.35cm}

\noindent{\bf Question (v).} It is indeed known that the {\em girth} of $G$ is infinite, that is, for 
suitable changes of the generating system, the length of the shortest nontrivial relation can be 
made arbitrarily large; see \cite{akhmedov-girthPL}. Is the girth of finitely-generated, non-solvable 
subgroups of $\mathrm{Diff}_+^{\omega}([0,1])$ infinite as well ?

\vspace{0.35cm}

Besides standard analytical methods used to deal with groups of diffeomorphisms of 
1-manifolds, the proof of Theorem A relies on a key argument to rule out laws in groups. 
This corresponds to a kind of ``finite ping-pong lemma'' that seems having been unexploited 
in full generality though widely known to the specialists (compare \cite[\S 4]{BS} and \cite{brin}). 
This strategy suits perfectly in case of relations of type (\ref{n-resilient-1}) and (\ref{n-resilient-2}) above. 

\vsp

\begin{rem} The group $\mathrm{Diff}^{\omega}_+([0,1])$ is residually solvable (just by truncating the Taylor 
series expansions of group elements at the origin). However, this is not enough to establish solvability in the case 
where some nontrivial law is satisfied, as shown by \cite[Theorem 8]{de-cornulier}. However, truncating series produces 
not only solvable but also algebraic (over the reals) quotiens.  An argument using Zariski closures may them be 
implemented to deduce Theorem A in a more algebraic (yet nontrivial) way. Again, it is hopeless to try to 
use a similar argument in lower regularity, hence we do not pursue this issue here and we leave 
the details to the reader.
\end{rem}

\vspace{0.4cm}


\noindent{\bf (b) An approach to Linnell's question and beyond.} The results of this subsection arose 
during several discussions with B. Deroin and V. Kleptsyn. They are almost direct consequences of 
\cite{random}, yet they deserve to be isolated in order to stress two major open questions in the subject.

As previously recalled, countable left-orderable groups are those that act on the real line by 
(orientation-preserving) homeomorphisms without a global fixed point. If the group is finitely-generated, 
such an action preserves a non-empty, closed minimal set, and is of one of the three types below:

\vspace{0.18cm}

\noindent I.- either it admits a discrete countable orbit or it is semiconjugate to the action of a dense group 
of translations; 
 
\vspace{0.18cm}

\noindent II.- it is not of type I and, up to a semiconjugacy if necessary, there is a 
homeomorphism of (semiconjugate) the real line with no fixed point that commutes 
with the (semiconjugate) action;

\vspace{0.18cm}

\noindent III.- it is of neither type I nor type II.

\vspace{0.18cm}

\noindent This classification is inspired from (the proof of) \cite[Th. 7.1]{random}. 
An almost direct consequence (details are given latter) is the next

\vspace{0.45cm} 

\noindent{\bf Theorem B.} {\em Let $\Gamma$ be a finitely-generated, left-orderable group having an action 
on the real line that is not of type III. Then either $\Gamma$ homomorphically surjects into $\mathbb{Z}$, or 
it contains a free subgroup in two generators.}

\vspace{0.45cm}

In particular, if there were no finitely-generated group all of whose faithful actions on the real line are of type 
III, then the answer to Linnell's question above would be affirmative. Besides this, there is another major 
question that would solve in the affirmative under the same assumption: No lattice in a higher-rank simple 
Lie group would embed into $\mathrm{Homeo}_+(\mathbb{R})$ (compare \cite{witte,lifschitz-witte}). Indeed, 
on the one hand, actions of type I lead to surjective homomorphisms into the reals just by taking translations 
numbers. On the other hand, actions of type II yield actions on the circle (viewed as the space of orbits of the 
commuting homeomorphism) that have no invariant probability measure (otherwise, we fall into type I). Both 
cases are impossible for higher rank lattices: the former contradicts Kazhdan's property (T), and the latter 
contradicts a theorem of Ghys \cite{ghys-reseaux}.  

\vsp

The discussion above shows the relevance of the next 

\vspace{0.35cm}

\noindent{\bf Question (vi).} Does there exist a finitely-generated, left-orderable group all of whose 
actions on the real line without global fixed points are of type III ?

\vspace{0.35cm}

Clearly, this is much related to the still open 

\vspace{0.35cm}

\noindent{\bf Question (vii).} Does there exist a finitely-generated, left-orderable, simple group ?

\vspace{0.45cm}


\noindent{\bf (c) Verbal properties of left-orders.} 
Our main result in this direction was obtained in collaboration with C. Rivas. Although 
not very surprising, it is by no means straighforward: The free group $F_2$ 
admits left-orders that satisfy no verbal property. Actually, this is the case 
of ``most'' left-orders on $F_2$, as we next explain. 

Recal that the set $\mathcal{LO}(\Gamma)$ of left-orders on a given left-orderable group $\Gamma$ 
carries a natural topology, where two orders are close if they coincide on a large finite set. This topology  
is metrizable if the group is countable. For instance, if $\Gamma$ is finitely generated, one way let 
$dist (\prec, \prec')$ be the inverse of the radius of the largest ball centered at the origin where 
$\prec$ and $\prec'$ coincide. The result below is stated (and proved) only for the free 
group in two generators, but it can be easily extended to the case of more generators. 

\vspace{0.45cm}

\noindent{\bf Theorem C.} {\em The set of left-orders on the free group $F_2$ satisfying 
no verbal property is a $G_{\delta}$-dense subset of the space of left-orders of $F_2$.} 

\vspace{0.45cm}

This result seems very different in nature to Theorems A and B. Nevertheless, as it will become  
clear along the proof, it relies on elementary combinatorial and dynamical arguments that are 
quite close to those involved in the proof of these two Theorems. Actually, for pedagogical 
reasons, we will prove our results in the reverse order.

\vspace{0.4cm}


\noindent{\bf Acknowledgments.} The author is strongly indebted to A.~Akhmedov, Y.~de~Cornulier, 
B.~Deroin, V.~Kleptsyn, and C.~Rivas, for quite useful and inspiring discussions, and for kindly allowing 
to include in this work several remarks and results that arose in these meetings. The author was funded 
by the Anillo Research Project 1103 DySyRF and the CNRS (UMR 8628, Univ. d'Orsay) via the ERC starting 
grant 257110 “RaWG”. He also acknowledges the IHP for the hospitality during the time this article was 
written, and would like to thank C.~Bleak, T.~Delzant, \'E.~Ghys, and R.~Strebel for their interest 
and useful references, as well as T.~Smirnova-Nagnibeda for the motivation to write this article. 


\section{Verbal Properties of Left-Orders}
\label{verbal-section}

\subsection{Left-orders on $F_2$ violating a prescribed verbal property} 

The construction of such orders is done using a very simple dynamical idea. To do this, recall that any 
group $\Gamma$ of orientation-preserving homeomorphisms of the real-line can be ordered by 
declaring $f \succ g$ if and only if $f(0) > g(0)$. This order is left-invariant yet not necessarily 
total. To make it total, we just need to consider more ``references points'' than the origin; 
see \cite{fourier} for details. 

Next, given a reduced word $W$ in two letters and carrying positive and negative exponents, we will 
construct two homeomorphisms of the real line $f, g$, both moving the origin to the right, such that 
the element $W(f,g)$ moves the origin to the left. Then the order on $\langle f,g \rangle$ defined 
above satisfies $f \succ id$, $g \succ id$, and $W(f,g) \prec id$. Via the (parhaps non-faithful) 
action of $\mathbb{F}_2 := \langle a,b \rangle$ given by $\phi(a) := f$ and $\phi(b) := g$, this 
induces a partial left-invariant order on $\mathbb{F}_2$, still denoted by $\prec$. If it is not 
total, we may consider a {\em convex extension}: if $\preceq'$ is any left-order on $F_2$, 
we define $\preceq^*$ by letting $c \succ^* id$ if and only if either $\phi(c) \succ id$, or  
$\phi(c)$ and $id$ are $\prec$-incomparable ({\em i.e.} $\phi(c)(0)=0$) and $c \succ' id$. 
Then $\preceq^*$ is a left-order on $F_2$ that satisfies $f \succ^* id$, $g \succ^* id$, 
and $W(f,g) \prec^* id$, as desired.

The construction of the desired action is done as follows. By interchanging $a$ and $b$ if 
necessary, we may assume that the word $W=W(a,b)$ writes in the form $W = W_1 a^{-n} W_2$, 
where $W_2$ is either empty or a product of positive powers of $a$ and $b$, the integer $n$ is positive, 
and $W_1$ is arbitrary. 
Let us consider two local homeomorphisms defined on a right neighborhood of the real line
such that $f(0) > 0$, $g(0) > 0$ and $W_2(f,g)(0) < f^n (0)$. This can be easily done by taking
$f(0) \gg g(0)$ and letting $g$ be almost flat on a very large right-neighborhood of the origin.
If $W_1$ is empty, just extend $f$ and $g$ into homeomorphisms of the real line. Otherwise,
write $W_1 = a^{n_k} b^{m_k} \ldots a^{n_2} b^{m_2} a^{n_1} b^{m_1}$, where all $m_i, n_i$ are
nonzero excepting perhaps $n_k$. The extension of $f$ and $g$ to a left-neighborhood of the
origin depends on the signs of the exponents $m_i, n_i$, and is done in a constructive manner.
Namely, first extend $f$ slightly so that $f^{-n} W_2 (f,g)(0)$ is defined and $f$ has a fixed point
$x_1$ to the left of the origin.  Then extend $g$ to a left-neighborhood of the origin so that
$g^{m_1} f^{-n} W_2(f,g)(0) < x_1$ and $g$ has a fixed point $y_1$ to the left of $x_1$. Notice
that $m_1 > 0$ forces $g$ to be right-topologically-attracting towards $y_1$ on an interval containg
$f^{-n} W_2(f,g)(0)$, whereas $m_1 < 0$ forces right topological repulsion. Next, extend $f$ to a left
neighborhood of $x_1$ so that $f^{n_1} g^{m_1} f^{-n} W_2(f,g)(0) < y_1$ and $f$ has a fixed point
$x_2$ to the left of $y_1$. Again, if $n_1 > 0$, this forces right-topological-attraction towards 
$x_2$, whereas $n_1 < 0$ implies right-topological-repulsion.

Continuing the procedure in this manner (see Figure 1 for an 
illustration), we get partially-defined  homeomorphisms $f,g$ for which
$$0 > f^{n_k} g^{m_k} \ldots f^{n_2} g^{m_2} f^{n_1} g^{m_1} f^{-n} W_2(f,g)(0) = W(f,g)(0).$$
Extending $f,g$ arbitrarily into homeomorphisms of the real line, we finally obtain the desired action.


\vspace{0.9cm}
\beginpicture
\label{picture verbal orderings}

\setcoordinatesystem units <1cm,1cm>

\plot 0 -2 0 2 / \plot  -4.3 0 2 0 /

\put{$f$} at 2.3 2

\plot 
0 1.55 
0.1 1.64 
0.2 1.69 
0.5 1.8 
1 1.89
1.5 1.96 
2 2 / 

\plot 
0 1.55
-0.1 1.23 
-0.2 0.35 
-0.3 0.03 
-0.4 -0.28  
-0.6 -0.64
-0.65 -0.65 / 

\plot 
-0.65 -0.65 
-0.7 -0.69 
-0.8 -0.78 
-0.9 -0.8
-2.2 -0.82
-2.3 -0.825
-2.4 -0.83
-2.5 -0.85 / 

\plot 
-2.5 -0.85 
-2.6 -1 
-2.7 -2.7 /
 
\plot 
-2.7 -2.7 
-2.71 -3 
-2.72 -3.2 
-2.75 -3.6 
-2.77 -3.78 
-2.78 -3.793  
-2.79 -3.796
-2.8 -3.9 / 

\plot 
-2.8 -3.9 
-3 -3.96
-3.4 -3.97 
-3.8 -3.98 
-4.3 -4 /

\put{$g$} at 2.4 0.4

\plot -0.05 0.2 
0 0.3 
2 0.4 / 

\plot -0.05 0.2 
-0.2 -1.3 
-0.25 -1.45
-0.3 -1.5 / 

\plot -0.3 -1.5 
-0.35 -1.55 
-0.4 -1.6 
-1.7 -1.7 
-3 -1.8 
-3.2 -1.86
-3.3 -1.9
-3.4 -2 / 

\plot 
-3.4 -2 
-3.45 -3
-3.47 -3.4 
-3.48 -3.48
-3.49 -3.49
-3.5 -3.5 / 

\plot 
-3.5 -3.5  
-3.8 -3.58
-4 -3.592 
-4.1 -3.595 
-4.2 -3.598
-4.3 -3.6 /

\put{{\tiny $\bullet$}} at  -0.65 -0.65 
\put{{\tiny $\bullet$}} at  -2.7 -2.7
\put{{\tiny $\bullet$}} at  -3.99 -3.99
\put{{\tiny $\bullet$}} at   -1.7 -1.7 
\put{{\tiny $\bullet$}} at  -3.5 -3.5


\setdots  \plot -4.3 -4.3 2 2 /

\small

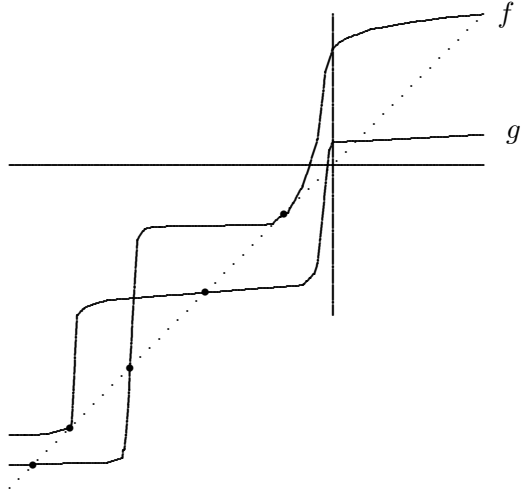
\put{Figure 1 \!: $W_1=a^{n_2}b^{m_2}a^{n_1}b^{m_1}$, where $m_1 >0, n_1 < 0, m_2 <0, n_2 >0$.} at -1 -5

\put{ } at -9.1 2.2
\endpicture


\vspace{0.5cm}

\subsection{Genericity of non-verbal orders} 

If $\Gamma$ is left-orderable, then it acts by conjugacy (equivalently, 
by right translations) on its space of left-orders. We denote by $\preceq_h$ the 
image of $\preceq$ under $h$ defined by letting $f \prec_h g$ if and only if $fh \prec gh$. 
Notice that given a word $W$, the subset of $W$-orders is preserved under the conjugacy 
action. Based on the work of McCleary \cite{Mc}, Clay \cite{clay} and, independently, Rivas 
\cite{rivas}, proved that $F_2$ carries left-orders whose orbits under the conjugacy action 
are dense. (This easily yields a new proof of the fact that the space of left-orders of $F_2$ 
is a Cantor set \cite{Mc,fourier}.) We next show that these orbits are made 
of orders satisfying no verbal property.

\vspace{0.1cm}

\begin{lem} {\em Every left-order on $F_2$ having a dense orbit 
under the conjugacy action satisfies no verbal property.}
\end{lem}

\noindent{\bf Proof.} Given a reduced word $W$ in two letters carrying positive and negative exponents, 
let $\preceq'$ be a left-order on $F_2$ that is not a $W$-order. Then there exist $\preceq'$-positive 
elements $f,g$ in $\Gamma$ such that $W(f,g)$ is $\preceq'$-negative. If $\preceq$ has a dense
orbit under the conjugacy action, then there exists $h \in \Gamma$ such that
$\preceq_h$ satisfies these three inequalities, namely, $f,g$ are
both $\preceq_h$-positive, whereas $W(f,g)$ is $\preceq_h$-negative. Hence,
$hfh^{-1}$, $hgh^{-1}$ are both $\preceq$-positive, whereas $h W(f,g) h^{-1} =
W(hfh^{-1},hgh^{-1})$ is negative. This shows that $\preceq$ does not satisfy the
$W$-verbal property.
$\hfill\square$

\vspace{0.4cm}

The proof of Theorem C can be now finished by a standard Baire type argument. Indeed, let us enumerate 
as $\{ W_1, W_2, \ldots \}$ all reduced words in two letters carrying positive and negative exponents. By the 
preceding Lemma, for each $W_i$, the set $\mathcal{LO}_{W_i} (F_2)$ of $W_i$-orders on $F_2$ has empty 
interior. By the definition of the topology on $\mathcal{LO} (F_2)$, this set is closed. Therefore, the 
complement of the union $\bigcup_i \mathcal{LO}_{W_i} (F_2)$ is a $G_{\delta}$-dense set. This 
complement corresponds to the set of left-orders satisfying no verbal property.  

\vspace{0.4cm}

\noindent{\bf Question (viii).} It is a nontrivial fact that the real-analytic homeomorphisms of the line 
given by $x \mapsto x+1$ and $x \mapsto x^3$ generate a free group \cite{glass-free}. By analyticity, 
a $G_{\delta}$-dense subset $S$ of points in the line have a free orbit under this action. Given a point 
$x \in S$, we may associate to it the left-order on $F_2$ defined by $f \succ g$ whenever $f(x) > g(x)$. 
Is the set of $x \!\in\! S$ for which the associate order satisfies no verbal property still a 
$G_{\delta}$-dense subset of $\mathbb{R}$ ?


\section{An Approach to Linnell's Question and Beyond}

As we have already suggested, the proof of Theorem B is somewhat tautological. Indeed, for actions of type I, 
a homomorphism into the reals is provided by the translation number. For actions of type II, we consider 
the action on the circle obtained as the quotient space of the commuting homeomorphism. This action 
cannot preserve a probability measure, otherwise the original action on the line would be of type I. 
Therefore, by the weak Tits alternative (see \cite{margulis,ghys}), the group contains free 
subgroups in two generators.

\vspace{0.2cm}

Because of the proof above, it becomes desirable a closer look at actions of type III. 
This is done by the next 

\vspace{0.1cm}

\begin{prop} \label{bluf}
{\em Let $\Gamma$ be a finitely-generated subgroup of $\mathrm{Homeo}_+(\mathbb{R})$. 
If $\Gamma$ acts with no global fixed point, then its action is of type III if and only if the 
following condition is satisfied: there exist $c < c'$ such that for all pairs of points 
$a <  b$ and $a' < b'$ in the real line such that $a < c < c' < b$, there is 
$g \in \Gamma$ such that $g(a) < a'$ and $g(b) > b'$.}
\end{prop}

\noindent{\bf Proof.} One direction is obvious: the presence of ``expanding elements'' prevents types I and II. To 
show the converse, assume the action of $\Gamma$ is of neither type I nor type II.  Let $K\!\subset\!\mathbb{R}$ 
be a minimal invariant closed set for the $\Gamma$-action (see \cite[Proposition 2.1.12]{book}). There are three 
cases:

\begin{itemize}

\item If $K$ is discrete, then it consists of a bi-infinite sequence of points ordered on the real line, which 
necessarily diverges in each direction (because $\Gamma$ admits no global fixed point). 
The action is hence of type I.

\item If $K$ is the whole line, then there are two subcases. 

-- If the action is free, then  by H\"older's theorem,  
it is conjugate to an action by translations (see \cite[\S 2.2.4]{book}), hence of type I. 

-- If it is non-free, then there must 
exist $g \in \Gamma$ and a point $x_0 \in \mathbb{R}$ which is fixed by $g$ and such that $g$ 
has no fixed point in an interval of type $[x_0-\varepsilon,x_0]$ or $[x_0,x_0+\varepsilon]$. 
In each case, the corresponding interval contracts into the point $x_0$ under iterates of $g$. 
As the action is minimal, for each $x \in \mathbb{R}$ we may consider the supremum 
$\varphi(x)$ of the $y \in (x,\infty)$ for which there exists a sequence $f_k \in \Gamma$ 
such that $f_k([x,y])$ converges to a single point in $\mathbb{R}$. Clearly, the function 
$\varphi \! : \mathbb{R}  \rightarrow \mathbb{R} \cup \{\infty\}$ is $\Gamma$-equivariant 
and nondecreasing.  If $\varphi (x) = \infty$ for some $x$, then this holds for all $x$. We claim 
that this allows performing the desired expansions. Indeed, as $\varphi(a') = \infty$, the interval 
$[a,b']$ can be contracted towards a point $y_0 \in \mathbb{R}$ by a sequence $f_k \in \Gamma$. 
As the action is minimal, there is $h \in \Gamma$ such that $h (y_0) \in (a,b)$. Then for a 
large-enough $k$, the element $g := (h f_k)^{-1}$ satisfies $g(a) < a'$ and $g(b) > b'$.

Assume next that $\varphi(x)$ belongs to $\mathbb{R}$ for all $x$. We claim that $\varphi$ 
is a homeomorphism of $\mathbb{R}$. Indeed, $\varphi$ has to be continuous, otherwise the 
interior of the complement of its image would be an open $\Gamma$-invariant set, which 
contradicts the minimality of the $\Gamma$-action. Surjectivity follows from that $\phi$ is 
nondecreasing and continuous, and injectivity is proved by contradiction: otherwise, the set 
of points having a neighborhood on which $\varphi$ is constant would be open and 
$\Gamma$-invariant,  thus contradicting minimality. 

Therefore, since $\varphi$ has no fixed point and is $\Gamma$-equivariant, we have that the action, 
if not of type I, is of type II.

\item Finally, if $K$ is not the whole line, then it is locally isomorphic to a Cantor set, 
so that the action is semiconjugate to that of a group acting minimally on $\mathbb{R}$. 
(Just by the classical trick of collapsing connected components of the complement into 
points.) We may hence apply the previous arguments, thus either showing that the 
action is semiconjugate to an action by translations or inducing elements in 
the original group that realize the desired expansions. We just need to take 
care in choosing both $c \!<\! c'$ lying in $K$ (or in different connected 
components of its complement). 
\end{itemize}\vspace{-0.832cm}
$\hfill\square$

\vspace{0.4cm}

\begin{ex} The free group $F_2$ admits actions of the three types above. Actions of type I come as 
dynamical realizations (in the sense of \cite[Teor. 6.8]{ghys-survey}) of bi-orders on $F_2$ (which do 
exist according to a classical result of Magnus; see e.g. \cite{DDRW}). Actions of type II are obtained just 
by lifting to the real line the generators of a Schotky group of circle homeomorphisms (this can be hence 
realized inside $\widetilde{\mathrm{PSL}}(2,\mathbb{R})$). Finally, actions of type III can be built ``by 
hand''. Actually, arguments as those of the previous section show that ``most actions'' of $F_2$ are 
of type III.
\end{ex}

\begin{ex} Actions of the Baumslag-Solitar group $BS(1,2) := \langle a,b \!: bab^{-1} = a^2 \rangle$ 
on the line (without fixed points) were classified up to semiconjugacy in \cite{rivas-conrad}. These 
all come from different inclusions in the affine group, except for four non-semiconjugate actions 
in which the element $b$ acts with no global fixed point. The former actions are of type III, 
whereas the latter are of type I. 
\end{ex}

\begin{ex} There are many examples of groups all of whose actions on the line 
are of type I. Concerning this, we may address the next

\vsp\vsp

\noindent{\bf Question (ix).} What are the left-orderable groups for all actions, the restriction to 
finitely-generated subgroups are of type I ?

\vsp\vsp

\noindent According to \cite{fourier}, this is equivalent to asking for the groups all of 
whose left-orders are Conradian. This includes left-orderable groups with no free subsemigroups 
(e.g. all torsion-free, nilpotent groups) as well as groups with finitely many left-orders \cite{KK}.
\end{ex}

An element $g$ of a left-orderable group $\Gamma$ is said to be {\em cofinal} if for 
every action on the real line with no global fixed point, $g$ fixes no point. Clearly, if a 
finitely-generated, left-orderable group $\Gamma$ has a cofinal, central element, 
then no action is of type III. 

\begin{ex} The lifting in $\widetilde{\mathrm{PSL}}(2,\mathbb{R})$ of the $(2,3,7)$ triangle 
group has the presentation 
$$G = \big\langle f,g,h \! : f^2 = g^3 = h^7 = fgh \big\rangle.$$
This example was introduced by Thurston in \cite{thurston}; it is the first example in the 
literature of a group of homeomorphisms of the real line with no nontrivial homomorphisms 
into the reals (hence no action of type I). The central element $\Delta := fgh$ is cofinal. (Indeed, 
if $\Delta = f^2 = g^3 = h^7$ has a fixed point, then this is fixed by $f,g,h$, hence by the 
whole group.) As a consequence, every action of $G$ is of type II. 
\end{ex}

\begin{ex} The center of the braid group is generated by the square of the so-called 
Garside element $\Delta_n$, which satisfies 
$$\Delta_n^2 = (\s_1 \s_2 \cdots \s_{n-1} )^n = (\s_1^2 \s_2 \cdots \s_{n-1})^{n-1}.$$
(Here, the $\sigma_i$'s are the canonical (Artin) generators of $B_n$.)  It was 
shown by Clay in \cite{cofinal} that $\Delta_n$ is cofinal in $B_n$. We do 
not know whether there exist type III actions of the commutator subgroups 
$[B_n,B_n]$ for $n \geq 5$ (these groups do not admit actions of type I since 
they admit no nontrivial homomorphism into the reals; the actions obtained by 
the Nielsen-Thurston method \cite{SW,NW} are of type II). Actually, we do not known 
any example of a finitely-generated, left-orderable group all of whose actions on 
the line with no global fixed point are of either type II or type III, and both types arise.
\end{ex}

\vsp\vsp

The classification into three types according to the dynamical properties 
can be made not only for actions on the line, but also for left-orders on 
a given left-orderable group. We will address this issue elsewhere. 


\section{Groups of Diffeomorphisms Satisfying Laws} 

\subsection{A ping-pong like lemma.} 

Below we give a ``finite version" of the classical ping-pong lemma of Klein (see \cite{harpe} 
for a discussion of the original version).

\vspace{0.1cm}

\begin{lem} \label{ping}
Let $\Gamma$ be a group acting by bijections on a set $X$. Assume that for every $k \in \mathbb{N}$, 
there exist elements $f,g$ in $\Gamma$ and nonempty subsets $A_1,\ldots A_k,B_1,\ldots,B_k$ of 
$X$ such that:\\

\noindent -- for all nonzero integers $n$, we have 
$f^n(A_i) \subset B_i$ for $i \!\in\! \{1,\ldots,k\}$, and 
$g^n(B_i) \subset A_{i+1}$ for $i \in \{1,\ldots,k-1\}$;\\  

\noindent -- the sets $A_1$ and $B_k$ are disjoint.\\

\noindent Then $\Gamma$ satisfies no nontrivial law.
\end{lem}

\noindent{\bf Proof.} As it is well-known, if a group satisfies a law, then it satisfies 
a law in two letters. (This easily follows from that the free group in two generators 
contain copies of free groups on arbitrarily many generators.) Let $W = W(a,b)$ be 
a word in two letters representing a group law. By conjugating $W$ by a 
power of $a$ if necessary, we may assume that it has the form 
$$W = a^{n_{k}} b^{m_{k-1}} a^{n_{k-1}} \cdots b^{m_1} a^{n_1},$$
where all exponents are nonzero. Consider the elements $f,g$ and the sets $A_i,B_i$ 
provided by the hypothesis for the integer $k$. We have 
\begin{eqnarray*}
W(f,g) (A_1) 
&=& f^{n_{k}} g^{m_{k-1}} f^{n_{k-1}} \cdots f^{n_2} g^{m_1} f^{n_1} (A_1)\\
&\subset& f^{n_{k}} g^{m_{k-1}} f^{n_{k-1}} \cdots f^{n_2} g^{m_1} (B_1)\\
&\subset& f^{n_{k}} g^{m_{k-1}} f^{n_{k-1}} \cdots f^{n_2} (A_2) \hspace{0.2cm}
\subset \hspace{0.08cm} \ldots \hspace{0.08cm} 
\subset \hspace{0.2cm} f^{n_k} (A_{k}) \hspace{0.12cm} 
\subset \hspace{0.12cm} B_{k}.
\end{eqnarray*} 
Since $A_1$ and $B_{k}$ are disjoint, this implies that $W(f,g)$ is a nontrivial 
element of $\Gamma$. Thus, $\Gamma$ does not satisfy the law given by $W$. 
$\hfill\square$

\begin{ex} \label{law-order} Assume that a group $\Gamma$ with a left-order $\preceq$ contains group 
elements satisfying (\ref{n-resilient-1}) and (\ref{n-resilient-2}). Given two group elements $f' \prec g'$, 
denote $[f',g']$ the set of group elements $h'$ satisfying $f' \preceq h' \preceq g'$. If we let 
$$A_i := \bigcup_{j=-i+1}^{i-1} \big[ h_{k+1+j}, h_{k+1+j}' \big], 
\qquad 
B_i := \bigcup_{j=-i}^{i-1} \big[ \bar{h}_{k+1+j}, \bar{h}_{k+1+j}' \big],$$
then conditions (\ref{n-resilient-1}) and (\ref{n-resilient-2}) translate into that the hypothesis of the 
preceding Lemma are satisfied. 
\end{ex}

\begin{ex} \label{ejemplo-F}
For the Thompson group $F$, given $k \!\in\! \mathbb{N}$, choose an element $g$ having 
exactly $2k+1$ fixed points inside $]0,1[$, all of them transversal. Denote (and number) these points so that 
$p_1 < \ldots < p_{2k+1}$. Let $f$ be another element having $2k$ fixed points $\{q_1, \ldots, q_{2k+1}\}$ 
in $]0,1[$, all of them transversal, so that $p_1 < q_1 < p_2 < q_2 \ldots < p_{2k+1} < q_{2k+1}$. (Such an  
$f$ can be taken as a conjugate of $g$, but this is irrelevant here.) Let $A$ and, respectively, $B$, be the 
union of (small-enough) disjoint neighborhoods of the points $p_i$ and $q_i$. By looking at all 
combinatorial possibilities, one can easily see that $A$ and $B$ contain subsets 
$A_1,\ldots,A_k,B_1,\ldots,B_k$ satisfying the conditions of Lemma \ref{ping} 
with respect to very large powers $f ^N$ and $g^N$. Indeed, for the set 
$A_1$ we take a small neighborhood of $p_{k+1}$, for $B_1$ the union of 
small neighborhoods of $q_{k}$ and $q_{k+1}$, for $A_2$ the union of 
small neighborhoods of $p_{k}$, $p_{k+1}$ and $p_{k+2}$, and so on. 
The reader should notice that this argument -which is nothing but a dynamical 
restatement of that of the preceding Example-- is not so far away from the proof 
in \cite[\S 5]{BS} that $F$ satisfies no law, despite the fact --also proved in \cite{BS}-- 
that it contains no free subgroup in two generators (see \cite{calegari} for a completely 
different proof of this).  
\end{ex}

\begin{rem} \label{remark-F}
In the example above, we do not really need that $f$ and $g$ have exactly $2k+1$ fixed points. 
What is essential is that they have at least $2k+1$ fixed points that are intertwined as above and 
that they admit no common fixed point in the interval $[p_1,q_{2k+1}]$. We leave the details of this 
to the reader.
\end{rem}


\subsection{A proof of Theorem A}

We start by recalling that a solvable group of real-analytic diffeomorphisms of the closed interval is 
necessarily metabelian, and its action is topologically conjugate to that of an affine group provided there 
is no global fixed point in the interior (see \cite{ghys}). The crucial step of the proof is given by the next

\vspace{0.2cm}

\begin{prop} \label{N-fixed-points}
{\em Let $\Gamma$ be a subgroup of $\mathrm{Diff}^{\omega}_+([0,1])$ for 
which there exists $N \geq 1$ such that every nontrivial group element has at most 
$N$ fixed points. Then $\Gamma$ is metabelian.}
\end{prop}

\vspace{0.2cm}

Let us mention that Akhmedov has recently shown a $C^{1+\alpha}$ version of this result; 
see \cite{akhmedov}. To keep this work reasonably self-contained, below we offer an 
elementary proof of the Proposition (in the real-analytic setting). This relies on next

\vspace{0.2cm}

\begin{lem}\label{nondiscrete}
{\em Let $\Gamma$ be a non-metabelian  
subgroup of $\mathrm{Diff}_+^{\omega}([0,1])$ having 
no global fixed point other than $0$ and $1$. Then for all points 
$0 < p_1 < \ldots < p_n < 1$ and all $\varepsilon > 0$, there is an 
element $f \in \Gamma$ that doesn't fix any of these points, though 
$| f(p_i) - p_i | < \varepsilon$ for all $i$.} 
\end{lem}

\noindent{\bf Proof.} Recall that non-metabelian groups of real-analytic 
diffeomorphisms of the interval contain local flows in their closure\footnote{Actually, 
this is rather elementary for groups acting on the interval: see Footnote 3.} 
\cite{Nakai,Rebelo}; in particular, they act minimaly on $(0,1)$ (see also 
\cite[Proposition 3.9]{ghys}). Hence, it suffices to show that every given interval can be 
mapped (expanded) into intervals whose endpoints are as close to 0 and 1 as desired. 
(Indeed, the conjugate flow will hence provide the desired elements for small translation 
parameters.) As we may certainly assume that $\Gamma$ is finitely generated, this 
amounts to saying that the action of $\Gamma$ on $(0,1)$ is of type III (see Proposition \ref{bluf}).

Suppose the action is of type I, that is, $\Gamma$ either preserves an infinite discrete set in $(0,1)$ 
or is semiconjugate to a group of translations. Then every element in the commutator subgroup 
$[\Gamma,\Gamma]$ fixes infinitely many points of $(0,1)$. By analyticity, these elements are 
necessarily trivial, hence $\Gamma$ is Abelian, which is a contradiction.

Assume the action is of type II. If $g \in \Gamma$ 
is nontrivial and has a fixed point, then the whole (infinite) orbit of this point under the homeomorphism 
that commutes with the action is made of points that are fixed by $g$, which contradicts analyticity. 
Therefore, the action of $\Gamma$ on $(0,1)$ is free, hence by H\"older's theorem the group is Abelian, 
which is again a contradiction. (Actually, this also contradicts the fact that the action was of type II, hence 
not of type I.)
$\hfill\square$

\vspace{0.38cm}

The proof of Proposition \ref{N-fixed-points} proceeds by induction on the maximum 
$N$ of the number $|max^{Fix}_{0}(\Gamma)|$ of fixed points in $(0,1)$ of a nontrivial 
element. If $N = 0$, then the action is free on the interior. By H\"older's theorem, 
the group $\Gamma$ is Abelian, and its action is conjugate to that of a group 
of translations. For $N = 1$, a theorem of Solodov establishes that $\Gamma$ 
must be metabelian, and its action is conjugate to that of an affine group 
\cite[\S 2.2.4]{book}. 
Assume that the claim holds for groups with nontrivial elements having at most 
$N-1$ fixed points, and let $\Gamma$ be a group with $|max^{Fix}_{0}(\Gamma)| = N$. To 
show that it is metabelian, we may suppose that it has no global fixed point other 
than $0$ and $1$. Assume for a contradiction that $\Gamma$ is non-metabelian. 
We claim that $\Gamma$ must contain: 

\vspace{0.18cm}

\noindent -- a nontrivial element $h$;

\vspace{0.18cm}

\noindent -- an element $f$ satisfying $f(x) > x$ for all $x \!<\! 1$ very close to $1$ 
and having an order of contact to the identity at the origin smaller than that of $h$ 
(in the sense that $f(x) / h(x)$ goes to $0$ as $x$ converges to the origin along 
points $x$ for which $h(x) \neq 0$);

\vspace{0.18cm}

\noindent -- an element $g$ with exactly $N$ fixed points, all of them topologically transversal, 
having an order of contact to the identity at the endpoint $1$ smaller than that of $f$ and such 
that $g(x) > x$ for all $x<1$ very close to 1. 

\vspace{0.18cm}

Indeed, the element $h$ comes from that $\Gamma$ is nontrivial. Let $r \geq 1$ be 
such that all elements in the $r^{th}$ commutator subgroup $\Gamma_r$ have 
contact order to the identity at the origin smaller than that of $h$. If $\Gamma$ is 
non-metabelian, then it is non-solvable; therefore,  $\Gamma_{r}$ cannot be Abelian, 
hence it acts non-freely. Let $f$ be a nontrivial element therein such that $f(x) > x$ 
for all $x<1$ very close to 1. Let $s \geq 1$ be such that the order of contact 
to the identity at 1 of all elements in $\Gamma_{r+s}$ is smaller than that 
of $f$. If $\Gamma$ is non-metabelian, then so is $\Gamma_{r+s}$, hence 
$|max^{Fix}_{0}(\Gamma_{r+s})| = N$ (otherwise, we could apply the induction hypothesis). 
Let $p_1,\ldots,p_N$ be the set of fixed points in $(0,1)$ of an element $\bar{g}$ 
therein having a maximal number of fixed points. Applying Lemma \ref{nondiscrete}, 
we get elements $g_k \in \Gamma_{r+s}$ that move but very little all these points. We 
claim that for a large-enough $k$, we can take our desired element $g$ as being one of 
$g_k \bar{g}$, $g_k^{-1} \bar{g}$. Indeed, on the one hand, for each $p_i$ at which $\bar{g}$ 
is topologically transversal, both $g_k \bar{g}$ and $g_k^{-1} \bar{g}$ have at least a fixed 
point close to it (if $k$ is large enough). On the other hand, for each $p_j$ at which $\bar{g}$ 
is tangent to the identity, one of $g_k \bar{g}$, $g_k^{-1} \bar{g}$ has at least two fixed points 
in a small neighborhood of $p_j$, and the other one no fixed point therein 
(again, this provided $k$ is large enough). Therefore,
$$\mathrm{Fix}_0 (g_k \bar{g}) + \mathrm{Fix}_0 (g_k^{-1} \bar{g}) \geq 2N,$$
where the inequality is strict if at least one of the fixed points detected above is 
not topologically transversal. This obviously implies our claim.
 
\vspace{0.2cm}

Having the elements $f,g,h$ at hand, we will next search for a contradiction. Let $a$ be the 
smallest fixed point of $h$ in $(0,1]$. Since $\Gamma$ has no global fixed point in $(0,1)$, 
up to conjugating $g$, we may assume that $\mathrm{Fix}_0 (g) \subset (0,a-\varepsilon)$ 
for a certain $\varepsilon > 0$. Changing $h$ by $h^{-1}$ if necessary, we may also assume 
that $h(x) < x$ for small $x$. Let $f_k := h^{-k} f h^{k}$. It is well-known that $f_k$ must 
uniformly converge to the identity\footnote{This argument may be also used to replace the 
previous use of the Nakai-Rebelo's theorem on the existence of local flows in the adherence 
of non-solvable groups.} on any interval $[0,a-\varepsilon]$ (see for instance 
\cite[Lemma 4.4]{quasim}). 

Now, since $f_k$ uniformly converges to the identity on any compact subinterval of $[0,a)$ and the interior 
fixed points of $g$ are transversal and contained in $(0,a)$, for a large-enough $k$ we have that the 
graph of $f_k$ must cross that of $g$ at least at $N$ points in $(0,a)$. Moreover, it has to cross that 
of $g$ near 1, as the contact order of $g$ at 1 is smaller than that of $f_k$. (See Figure~2.) 
Therefore, the element $f_k^{-1} g$ (is nontrivial and) has at least $N+1$ fixed points, 
which is a contradiction. This closes the proof of Proposition \ref{N-fixed-points}.

\vspace{0.5cm}

\beginpicture

\setcoordinatesystem units <0.9cm,0.9cm>

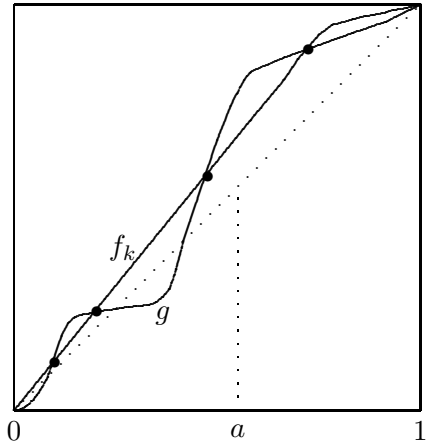
\put{Figure 2: The graphs of $f_k$ and $g$ cross $N=3$ times in $(0,a)$ and once more close to 1.} 
at 3 -0.94
\put{} at -5.8 0 

\putrule from 0 0 to 6 0
\putrule from 0 6 to 6 6
\putrule from 0 0 to 0 6   
\putrule from 6 0 to 6 6



\put {$0$}  at 0     -0.3 
\put {$a$}  at 3.3     -0.3
\put {$1$} at 6     -0.3

\put {$g$} at 2.2  1.4
\put {$f_k$} at 1.6  2.4 


\plot 
0 0
0.07071 0.01 
0.1 0.02 
0.15811 0.05 
0.2 0.08 
0.2236 0.1 
0.25495 0.13 
0.28482  0.16    
0.31622 0.2
0.33911 0.23 
0.36055 0.26 
0.38729 0.3 
0.41833 0.35 
0.44721 0.4 
0.47434 0.45 
0.5 0.5 /

\plot 
0.5 0.5 
0.7 1.1
0.8 1.29 
0.9 1.38
1 1.418
1.1  1.428
1.2 1.456 
1.3 1.47
1.5  1.5 /

\plot 
1.5 1.5 
1.6 1.52
1.7 1.53
1.8 1.54
1.9 1.55
2  1.56
2.1 1.6
2.2 1.68
2.3 1.8
2.4 2.1
2.5 2.5 /

\plot 
2.5 2.5
3 4
3.2 4.5 
3.3 4.7
3.4 4.89
3.5 5.01
4 5.22
4.5 5.4
5 5.58 
5.5 5.75
6 6 /

\plot 
0 0 
3.5 4.3 
4 4.9 
4.2 5.2 
4.5 5.55 
4.7 5.7 
5 5.78
5.2 5.835
5.3 5.845
5.5 5.885
5.6 5.92
5.8 5.95
6 6 /

\setdots
\putrule from 3.3 0 to 3.3 3.3
\plot 0 0 
6 6 / 

\begin{small}
\put{$\bullet$} at 0.47 0.7
\put{$\bullet$} at 1.09 1.445 
\put{$\bullet$} at 2.73 3.45 
\put{$\bullet$} at 4.22 5.33 

\end{small}

\endpicture


\vspace{0.5cm}

\noindent{\bf Proof of Theorem A.} 
Let $\Gamma$ be a non-metabelian subgroup of $\mathrm{Diff}_+^{\omega}([0,1])$. 
By Proposition \ref{N-fixed-points}, for every $k \!\geq\! 1$, the group $\Gamma$ contains 
an element $f$ having at least $4k+2$ fixed points at the interior; actually, according to the 
proof, these points can be supposed to be transversal. Using Lemma \ref{nondiscrete}, we 
may slightly perturb $f$ (just by left composition with an element that moves very little these 
points)  into a certain $g \in \Gamma$ still having $4k+2$ transversal fixed points but all 
different from those of $f$. After such a perturbation, we may easily find a sequence of at 
least $2k+1$ intertwined fixed points $p_1<q_1<p_2<q_2 \ldots < p_{2k+1} < q_{2k+1}$ 
\hspace{0.1cm} of $f$ and $g$, respectively. Using Lemma \ref{ping}, a modification 
of the argument of Example \ref{ejemplo-F} (see Remark \ref{remark-F}) then 
shows that $\Gamma$ satisfies no law.


\begin{footnotesize}

\vspace{0.2cm}

Andr\'es Navas\\

Univ. de Santiago de Chile, Alameda 3363, Santiago, Chile\\

email: andres.navas@usach.cl

\end{footnotesize}

\end{document}